\documentclass[11pt,leqno]{article}
\usepackage{amssymb,amsmath, amsthm}
\theoremstyle{definition}

\theoremstyle{definition}

\theoremstyle{definition}

\theoremstyle{remark}

\usepackage{amssymb, palatino}
\usepackage{graphicx}

\begin{document}
\title{ Minimal volume and simplicial norm of visibility n-manifolds and  compact 3-manifolds}
\author{
Jianguo Cao\thanks{Supported partially by NSF grant DMS-0706513 and
visiting Changjiang Chair Professorship at Nanjing University,
China.} \quad and \quad Xiaoyang Chen }
\date{}
\maketitle
\begin{abstract}
{\small In this survey paper, we shall derive the following result.

\noindent \textbf{Theorem A.} {\itshape Let $M^n$ denote a closed
Riemannian manifold with nonpositive sectional curvature and let
$\tilde M^n$ be  the universal cover of  $M^n$ with the lifted
metric. Suppose that the universal cover $\tilde M^n$ contains no
totally geodesic embedded Euclidean plane $\mathbb{R}^2$ (i.e.,
$M^n$ is a visibility manifold ).  Then Gromov's simplicial volume
$\| M^n \|$ is non-zero. Consequently, $M^n$ is non-collapsible
while keeping Ricci curvature bounded from below. More precisely, if
$Ric_g \ge -(n-1)$, then $vol(M^n, g) \ge \frac{1}{(n-1)^n n!} \|
M^n \| > 0. $}

Among other things, we also outline a proof for the following direct
consequence  of Perelman's recent work on 3-manifolds.

\noindent \textbf{Theorem B.} (Perelman) {\itshape  Let $M^3$ be a
closed a-spherical 3-manifold ($K(\pi, 1)$-space) with the
fundamental group $\Gamma$. Suppose that $\Gamma$ contains no
subgroups isomorphic to $\mathbb{Z}\oplus \mathbb{Z}$. Then $M^3$
is diffeomorphic to a compact quotient of real hyperbolic space
$\mathbb{H}^3$, i.e., $M^3 \equiv \mathbb{H}^3/\Gamma$.
Consequently, $MinVol(M^3) \ge \frac{1}{24}\| M^3 \| > 0$. }

Minimal volume and simplicial norm of all other compact
3-manifolds without boundary and {\it singular} spaces will also be discussed.}

\end{abstract}


\section*{1. Introduction.}
\noindent
 In 1993, Professor S. T. Yau proposed studying the Martin boundary, the Gromov-norm and other semi-hyperbolic properties of manifolds with Ballmann rank one
 (cf. Problem \#47 and \#97 of [Y93]).
 A nonpositively curved manifold $M^n$ is called to be of  Ballmann rank one, if there is a geodesic $\sigma: \mathbb{R} \to M^n$ which does not admit any
 non-zero orthogonal parallel Jacobi field (cf. [Ba82]). Various progress has been made in the study of rank-one manifolds, for example, see
 [BaL94], [BCG96], [Ca00],
 [CCR01], [CCR04], [CFL07] and [LS06].  Manifolds of Ballmann rank-one can be divided into two sub-classes: collapsible manifolds and non-collapsing manifolds.
 For collapsing manifolds of nonpositive sectional curvature, the first author, Jeff Cheeger and Xiaochun Rong indeed showed that
 {\it ``if a closed Riemannian manifolds with nonpositive sectional curvature
 $-1 \le sec_{M^n} \le 0$ and if the volume is sufficiently small $vol(M^n) \le \varepsilon_n$, then $M^n$ must be a generalized graph-manifold and hence
 minimal volume $MinVol(M^n)$ and Gromv-norm $\| [ M^n]\|$ of $M^n$ both vanish"}, (cf. [CCR01], [CCR04]). In this paper, we consider
 non-collapsible manifolds of Ballmann rank-one. In particular, we consider {\it visibility manifolds} in the sense of Eberlein and
 O'Neill. Our results of this paper will be complementing to the results of [CCR01-04] on collapsing manifolds with nonpositive curvature.

 In 1973 Eberlein and O'Neill introduced the so called "visibility
  manifolds".

  \bigskip

\noindent
\textbf{Definition 1.1.} ([EO73])
{\itshape
Let $M^n$ denote a closed Riemannian manifold with nonpositive
sectional curvature and let $\tilde M^n$ be  the universal cover
of  $M^n$ with the lifted metric. If the universal cover $\tilde
M^n$ contains no totally geodesic embedded Euclidean plane
$\mathbb{R}^2$, then $M^n$ is called a visibility manifold.}

\bigskip
   In this paper, we always assume that all spaces have dimensions $n\geq2$. It has been conjectured by various authors  that the following assertion might be true.

   \vspace{5mm}
   \noindent
   \textbf{ Conjecture 1.2.} {\itshape Any closed visibility manifold $M^n$ must  admit a metric $g$ of negative sectional curvature $sec_g \le -1$.}
\vspace{5mm}

  There is partial evidence to support this conjecture. For example,
the Preissman theorem states that if a closed manifold
  $M^n$ admits a metric of negative sectional curvature, then its fundamental group $\pi_1(M^n)$ has no subgroup isomorphic to
$\mathbb{Z}\times\mathbb{Z}$, (cf. [Y71]). In fact, it follows from Yau's thesis [Y71] that
  there are no $\mathbb{Z}\times\mathbb{Z}$ subgroups contained
  in the fundamental group of a visibility manifold.

  In this paper, we will provide another piece of evidence for the
  conjecture.
  It is known that if a closed manifold $M^n$ admits a metric of negative sectional curvature, then $M^n$ is non-collapsible, (cf. [Gr82] and [IY82]).
  We shall also show that any closed visibility manifold is {\it non-collapsible}. In fact, Gromov introduced the notion of minimal volume by setting:
$$
  MinVol(M^n) = \inf \{ vol_n(M^n, g) \quad | -1 \le sec_g \le 1 \}. \eqno(1.1)
  $$

  More precisely, we shall show that any closed visibility manifold has non-zero Gromov-norm,
  non-zero minimal volume and hence is non-collapsible.
   In a seminal  paper ([Gr82]), Gromov introduced
  the notion of simplicial volume of a closed oriented manifold.
  In the same paper,
the question was raised which manifolds have non-zero
  simplicial volume. For a closed manifold of negative curvature, Inoue and Yano ([IY82])
  showed that the simpilicial volume must be nonzero. It is a natural
  question whether the simpilicial volume of visibility manifolds is
  nonzero. In this paper, we show that it turns out to be  true. More precisely,
  we prove the following theorem:

\medskip

\noindent
\textbf{Main Theorem.}
{\itshape
Let $M^n$ be a closed manifold of nonpositive curvature. Suppose
that  $M^n$ is a visibility manifold. Then Gromov's simplicial volume $\| M^n \|$ is non-zero.
Consequently, $M^n$ is non-collapsible while keeping Ricci curvature bounded from below. More precisely, if
$Ric_g \ge -(n-1)$, then
$$
vol(M^n, g) \ge \frac{1}{(n-1)^n n!} \| M^n \|  > 0. \eqno(1.2)
$$}
\medskip

We are very grateful to Igor Belegradek for bringing the work of
Bridson [Bri95] and Yamaguchi [Yama97] on singular spaces to our
attention. The definition of curvatures for possibly singular spaces
can be found in many graduate textbooks, (e.g. [BuBI01]).

\bigskip

\noindent \textbf{Main Corollary.} {\itshape Let $X$ be a
simply-connected and complete metric space of non-positive
curvature. Suppose that there does not exist an isometric embedding
of the Euclidean plane into $X$; $X$ has no boundary and suppose
that $X$ has a co-compact lattice $\Gamma$ such that $X/ \Gamma$ has
co-homological dimension $n$ and $H_n(X/\Gamma, \mathbb Z) \neq 0$.
Then the simplicial norm of $X/\Gamma$ is positive. If, in addition,
$X$ has curvature $\ge -1$, then
$$
vol(X/\Gamma) \ge \frac{1}{(n-1)^n n!} \| X/\Gamma \|  > 0.
\eqno(1.3)
$$
}

For spaces of dimension $ 3$, the work of Perelman provides a more
refined result to support Conjecture 1.2
  above as well.

\bigskip

\noindent \textbf{Theorem 1.3.} (due to Perelman [Per02-03],
[CZ06], [KL06], [GT07-08]) {\itshape  Let $M^3$ be a closed
aspherical 3-manifold ($K(\pi, 1)$-space) with the fundamental
group $\Gamma$. Suppose that $\Gamma$ contains no subgroups
isomorphic to $\mathbb{Z}\oplus \mathbb{Z}$. Then $M^3$ is
diffeomorphic to a compact quotient of real hyperbolic space
$\mathbb{H}^3$, i.e., $M^3 \equiv \mathbb{H}^3/\Gamma$.
Consequently, any compact 3-dimensional visibility manifold $M^3$
admits a metric of constant negative sectional curvature. }

\bigskip

Perelman posted two important preprints on Ricci flows on compact 3-manifolds with surgery online ([Per02], [Per03]), in order to solve
Thurston's conjecture on Geometrization of 3-dimensional manifolds. Thurston's Geometrization Conjectures states that ``for any closed, oriented
and  connected 3-manifold $M^3$, there is a decomposition $[M^3 -\bigcup \Sigma^2_j] = N^3_1 \cup N^3_2 ... \cup N^3_{m}$ such that each $N^3_i$
admits a locally homogeneous metric with possible incompressible boundaries $\Sigma^2_j$, where $\Sigma^2_j$ is homeomorphic to a quotient of a
2-sphere or a 2-torus". There are exactly 8 homogeneous spaces in dimension 3. The list of 3-dimensional homogeneous spaces includes 8
geometries: $\mathbb{R}^3$, $\mathbb{H}^3$, $\mathbf{S}^3$, $\mathbb{H}^2 \times \mathbb{R}$, $\mathbf{S}^2 \times \mathbb{R}$, $\tilde{SL}(2,
\mathbb{R})$, $Nil$ and $Sol$.  Several teams of outstanding mathematicians generously made their efforts to fill in the detailed proof of
Perelman's densely written arguments in [Per02] and [Per03]. Among them, we should mention important contributions in the clarification of
Perelman's work by
 Cao-Zhou [CZ06], Kleiner-Lott [KL06] and Morgan-Tian [MT07]-[MT08]. The Ricci flow with surgery part of Perelman's work has been well-understood.
 For example,
 see Ben Chow and his co-authors' four different books on Ricci flows (e.g. [CLN06]). The metric
 geometric part of Perelman's work (Perelman's collapsing theorem for 3-manifolds) has been studied
 by four  teams, including Besson's team [BBBMP07],
 Cao-Ge [CG09], Morgan-Tian [MT08] and Shioya-Yamaguchi [ShiY05].
 The papers [BBBMP07] and [MT08] did not use Perelman's stability theorem (cf. [Ka07]), while
 the paper of the first author and Ge [CG09] used the critical point theory of distance functions. Shioya and
  Yamaguchi used Perelman's stability theorem in their papers [ShiY00] and
 [ShiY05] before the publication of Vitali Kapovitch's proof [Ka07]. Over all, no major errors
 have been
 found in Perelman's papers [Per02]-[Per03]. Perelman's work is
 very solid and has become increasingly appreciated
by the geometry-topology community.
We will explain how Perelman's recent work ([Per02] and [Per03]) on Geometrization of 3-manifolds implies the above Theorem 1.3 in Section 4
below. Gromov's norm of other compact 3-manifolds will also be discussed in Section 4.

There are several ways to obtain the positivity of simplicial
volume. Following a suggestion of Thurston, it suffices to get a uniform upper bound
of the volumes of suitably defined straighted simplices. Gromov
[Gr82] outlined, and Inoue-Yano [InY82] showed in detail that  the
simplicial volume of closed negative curved manifolds is nonzero.
The recent work of Jean-Fran\c{c}ois Lafont and Benjamin Schmidt
verified the positivity of simplicial volume of a closed locally
symmetric space of non-compact type ([LS06]) by using a similar
ideas. F. Ledrappier and S. Lim recently studied volume entropy of
hyperbolic buildings, cf. [LeL08].

Our proof  takes a different approach. We will use Gromov's bounded cohomology method,  Gromov-hyperbolic group theory [Gr87] and
 a theorem of Miniyev [M01] to verify our Theorem. According to Gromov,
 to verify that the  simplicial
volume of $M^n$ is positive, it suffices to show that $M^n$ has nonzero $n$-dimensional bounded cohomology.
 A theorem of Igor Miniyev asserts that,
if the fundamental group $\pi_1(M^n)$ of $M^n$ is Gromov-hyperbolic and if $M^n$ is a closed oriented aspherical manifold ( a $K(\pi,
1)$-space), then $M^n$ has nonzero $n$-dimensional bounded cohomology. Finally, we will appeal to an earlier result of the first author, which
states that, for a closed manifold $M^n$ with nonpositive sectional curvature, $M^n$ is a visibility manifold if and only if its fundamental
group $\pi_1(M^n)$ is Gromov-hyperbolic, cf [Ca95] and [Ca00].  This  completes
the outline of  proof of our main theorem.

   The organization of this short article goes as follows. In Section 2, we review basic properties of similicial volume and bounded cohomology needed for our
   proof. In Section 3, we discuss the relations among visibility manifolds and Gromov-hyperbolic
    groups. The proof of the Main Theorem will be
   completed in Section 4.

\section*{2. Similicial volume and bounded cohomology.}

\noindent
  We begin with the definition of simplicial volume.
\bigskip

\noindent
\textbf{Definition 2.1.}
{\itshape
Let $M^n$ be a $n$-dimensional closed manifold, $C^0(\Delta^k, M^n)$
be the set of singular $k$-simplices and $C_k(M^n,\mathbb{R})$ be
the set of singular $k$-chains with real coefficient. For any
$c=\displaystyle \sum^j_{i=1}r_if_i$ with each $r_i\in\mathbb{R}$
and $f_i\in C^0(\Delta^k,M^n)$ be a singular real chain, the $l^1$-norm of c is defined by $\|c\|_1=\displaystyle\sum_{i=1}^j|r_i|$.
The $l^1$-norm of a real singular homology class $[\alpha]$ is
defined by
$$\|[\alpha]\|_1=inf\{\|c\|_1:\partial(c)=0,[c]=[\alpha]\}. \eqno(2.1)
$$ }

 To study the volume collapsing property of a given manifold, we concentrate on the norm of the top-dimensional
cohomology class.

\bigskip
\noindent
\textbf{Definition 2.2.} (Gromov's simplicial norm)
{\itshape Let $M^n$ be an oriented closed connected $n$-dimensional manifold.
The simplicial volume of $M^n$ is defined as
$\|M^n\|=\|i([M^n])\|_1$, where $i$:
$H_n(M,\mathbb{Z})\longrightarrow H_n(M,\mathbb{R})$ is the change
of coefficient homomophism, and $[M^n]$ is the fundamental class
arising from the orientation of $M^n$.}
\bigskip

   The following fundamental theorem of Gromov indicates that there are deep relations between minimal volume (defined by equation (1.1) above)
   and Gromov's simplicial norm.

\bigskip

\noindent \textbf{Theorem 2.3.} (Gromov's estimate [Gr82, page 220])
{\itshape Let $(M^n, g)$ be a complete $n$-dimensional Riemannian
manifold with Ricci curvature bounded from below by $Ric_g \ge
-(n-1)$. Then the volume of $(M^n, g)$ is bounded below by
$\frac{1}{(n-1)^n n!} \| M^n \|$. Consequently, one has
$$
  MinVol(M^n) \ge \frac{1}{(n-1)^n n!} \| M^n \|. \eqno(2.2)
$$
}
 Erwann Aubry further shown a similar estimate holds if the Ricci curvature of $M^n$ is
$L^p$ bounded from below by $-(n-1)$. We state his estimate as
follows:
\bigskip

\noindent \textbf{Theorem 2.4.} (Erwann Aubry's estimate [Au08])
{\itshape Let $M^n$ be a complete $n$-dimensional Riemannian
manifold. If for any $\epsilon>0$, $\exists$ a constant
$C(p,n,D,\epsilon)>0$ such that $(M^n, g)$ satisfies $D^2
\|\rho\|_p\leq C(p,n,D,\epsilon)$, then
$$ vol(M^n,
g) \ge \frac{1}{(n-1)^n n!} \| M^n \| . \eqno(2.3) $$} Here $D$
denotes the diameter of $M^n$ and $p > \frac{n}{2}$,
$\|\rho\|_p=(\frac{1}{Vol M^n}\int_M \rho^p)^\frac{1}{p}$,
$\rho=(\hat{Ric}+(n-1))_-$, $\hat{Ric}(x)$ is the least eigenvalue
of the Ricci tensor $Ric$ at $x$ and $f_- =max(-f,0)$.

When $M^n$ is diffeomorphic to a locally symmetric space of negative sectional curvature, Besson-Courtois-Gallot improved (2.2) by considering
metrics of scalar curvature bounded from below. We thank Professor D. Kotschick [Ko04] for bringing [BCG96] to our attention again.

\bigskip

\noindent \textbf{Theorem 2.5.} (Besson-Courtois-Gallot [BCG96])
{\itshape  Let $M^n$ be diffeomorphic to a compact locally
symmetric space of negative sectional curvature and $g_{can}$ be
the canonical metric of scalar curvature $ -n(n-1)$. Suppose that
$g$ is a smooth metric of scalar curvature $\ge -n(n-1)$ on $M^n$.
Then
$$
  Vol(M^n, g) \ge Vol(M^n, g_{can}) = c_n \|M^n\| > 0. \eqno(2.4)
$$
 }

We need to borrow another brilliant idea of Gromov to calculate the simplicial norm via the  dual-space method.
It is clear that the $\ell^\infty$-space is dual to the $\ell^1$-space.
Gromov introduced the so-called bounded cohomology as follows.


\bigskip

\noindent \textbf{Definition 2.6.} {\itshape Let
$C^k(M^n,\mathbb{R})$ be the set of singular $k$-cochains with real
coefficient. For any $c\in C^k(M^n,\mathbb{R})$, set the $l^\infty$
norm of $c$ by $|c|_\infty=sup\{|c(\sigma)|:\sigma \in
C_k(M^n,\mathbb{R})\}$. The set of bounded $k$-cochains with real
coefficient $C^k_b(M^n,\mathbb{R})$ is defined by
$$
C^k_b(M^n,\mathbb{R})=\{c: c \in C^k(M^n,\mathbb{R}),|c|_\infty <
\infty\}. \eqno(2.5)
$$
Assume that $\delta$ is the co-chain operator from
$C^k(M,\mathbb{R})$ to $C^{k+1}(M,\mathbb{R})$. It can be easily
checked that $\delta (C^k_b(M,\mathbb{R}))\subseteq
C^{k+1}_b(M,\mathbb{R})$ and we define the $k$-th bounded cohomology
group of $M^n$ as}
$$
H^k_b(M,\mathbb{R})=Ker(\delta (C^k_b(M,\mathbb{R}))/Im(\delta
(C^{k-1}_b(M,\mathbb{R})). \eqno(2.6)
$$

\bigskip

\noindent \textbf{Remark 2.7.}
If $G$ is an arbitrary group, the bounded cohomology group of $G$
can be defined as the corresponding bounded cohomology group of
the Eilenberg-Maclane space $K(G,1)$. For a further discussion of
the bounded cohomology theory of groups, see [Iv87] and [Bro94].

\bigskip

  As shown by Gromov, the theory of bounded cohomology sheds new
light on the notion of simplicial volume. In fact, it gives the
cohomological definition of simplicial volume. More precisely, we have
the following proposition.


\bigskip

\noindent \textbf{Proposition 2.8.} (Gromov [Gr82]) {\itshape Let
$M^n$ be  a closed oriented $n$-dimensional manifold. Then
 $ \|M\|^{-1}=inf\{\|\beta\|_\infty: \beta\in H^k_b(M,\mathbb{R}),
[\beta,[M]]=1\}$. Consequently,  $\|M\|$ is nonzero if and only if
there exists a bounded $\beta\in H^n_b(M,\mathbb{R})$ which does not
vanish on $[M]$.}

\bigskip

For the proof of the above proposition, see [Gr82].

In the next section, we set the stage to show that {\it ``if
$M^n$ is a closed visibility manifold then $\pi_1(M^n) $ must be Gromov-hyperbolic and hence $H^n_b(M^n,\mathbb{R}) \neq 0$"}.

\section*{3. Visibility manifolds and Gromov-hyperbolic spaces.}

\noindent
 There are at least three equivalent definitions of visibility manifolds. We will also use the following one of equivalent definitions for
 visibility manifolds.

\bigskip

\noindent \textbf{Definition 3.1.} ([EO73]) {\itshape (1) A simply-connected manifold $\tilde M^n$ of nonpositive curvature is said to be a
visibility manifold if for each point $p\in \tilde M^n$ and $\varepsilon>0$, there exists a constant $R(p,\varepsilon)>0$ such that if $\sigma:
[a,b]\rightarrow \tilde M^n$ is a geodesic segment satisfying the condition $d(p,\sigma)\geq R(p,\varepsilon)$, then
$\measuredangle_p(\sigma(a),\sigma(b))\leq\varepsilon$, where $\measuredangle_p(\sigma(a),\sigma(b))$ denotes the angle based at
$p$ between
$\sigma(a)$ and $\sigma(b)$. $\tilde M^n$ is called a "uniform visibility manifold" if the constant $R(p,\varepsilon)$ may be chosen to be
independent of $p$.

(2)  A closed manifold $M^n$  of nonpositive curvature is said to be a visibility (uniform visibility) manifold if its universal covering does.
}

\bigskip

 Examples of uniform visibility manifolds include all complete Riemannian manifolds with sectional curvature $\le -a^2 < 0$.

\smallskip

\noindent
\textbf{Example 3.2.} {\it (Negative curvature implies visibility). } A simply-connected and complete Riemannian manifold $\tilde M^n$
with strictly negative sectional
curvature $sec_{\tilde M^n} \leq -1 $ is a uniform visibility manifold. Since $\tilde M^n$ does not contain any totally geodesic
flat $\mathbb{R}^2$,  by Definition 1.1 in Section 1, we see that
$\tilde M^n$ is a visibility manifold. By an explicit calculation below,  we can further show that $\tilde M^n$ is a uniform visibility
manifold described in Definition 3.1. In fact, when $sec_{\tilde M^n} \leq -1 $, we can choose that
$$ R(p,\varepsilon) = R( \varepsilon)
= |\cosh^{-1} (\frac{ \pi}{\varepsilon} + 1)|.  \eqno(3.1)$$

To see this, we use two elementary facts. For $p\in \tilde M^n$ and a geodesic segment
$\sigma: [a,b]\rightarrow \tilde M^n$, we let $\Delta_{p, \sigma}$ be the cone with the apex $p$ over $\sigma([a, b])$. It is easy to
check that $\Delta_{p, \sigma}$ is a ruled sub-manifold of $\tilde M^n$. Let $\hat{g}_\Delta$ be the induced metric on the solid triangle
$\Delta_{p, \sigma}$.
There is a theorem of Gauss on the relation between intrinsic curvature and extrinsic curvature. Let $II_{\Delta}$ be the second fundamental form of
$\Delta$ in $\tilde M^n$. Suppose that $X$ is the unit ruled direction of $\Delta$ and $Y$ is its orthonormal complement of $X$ in $T_x(\Delta)$.  Then
a calculation shows that
$$
sec_{(\Delta,\hat g) } = \sec_{\tilde M^n} + \langle II_{\Delta}(X, X), II_{\Delta}(Y, Y) \rangle - \| \langle II_{\Delta}(X, Y)\|^2 \le -1, \eqno(3.2)
$$
for all $x \in \Delta$.

Using (3.2) and applying Gauss-Bonnet formula to the geodesic triangle $(\Delta,\hat g)$, we have the area estimate of $(\Delta,\hat g)$:
$$
Area(\Delta,\hat g) \le \pi. \eqno(3.3)
$$

Let $B(p, r)$ be the metric ball in $\tilde M^n$ centered at $p$ with radius $r$, and let $\varepsilon = \angle_p(\sigma(a), \sigma(b))$ be the
visibility angle
of $\sigma([a, b])$ viewed from $p$. Let $R$ be the largest  radius of inscribed fan in
$\Delta$ given by
$$
  R = \inf_{a \le t \le b} \{ \| Exp_p^{-1}(\sigma(t))\|  \}, \eqno(3.4)
$$
where $Exp_p$ is the exponential map of $\tilde M^n$ at $p$.

 Using area comparison theorem for $\sec \le -1$, we have
$$
[\cosh R -1] \varepsilon \le Area[B(p, R) \cap \Delta_{p, \sigma}]\le  \pi. \eqno(3.5)
$$
The formula (3.1) now follows from (3.4)-(3.5). Hence, when $\tilde M^n$ has strictly negative sectional curvature
$sec_{\tilde M^n} \leq -1 $, such a manifold $\tilde M^n$ must be
a uniform visibility manifold with the visibility function
$$   R(p,\varepsilon)=  R(\varepsilon)
= |\cosh^{-1} (\frac{ \pi}{\varepsilon} + 1)|.
$$
This completes the proof of the assertion that ``strictly negative curvature implies uniform visibility". \qed
\smallskip

  It is a very interesting fact that a uniform visibility manifold turns out
to be a Gromov hyperbolic space, which is very important in the
study of large scale geometry. Roughly speaking, Gromov hyperbolic
space characterizes coarse features of hyperbolic plane, and can be
extended to the category of metric space.
We need to recall some basic properties of Gromov-hyperbolic theory, which are needed in our paper.
 For a detailed exposition of Gromov-hyperbolic groups and Gromov-hyperbolic spaces, see [Sho91].

\bigskip

\noindent
\textbf{Definition 3.3.}
{\itshape
Let $(X,d)$ be a metric space with underlying set $X$ and metric
$d$, and let $x,y\in X, \gamma:[0,1]\rightarrow X$ be a curve in $X$
connecting $x$ and $y$, define the length of $\gamma$ by
$$
\ell(\gamma)=sup\{\displaystyle\sum_{i=0}^n d(x_i,x_{i+1}),
x=x_0<x_1<\cdots<x_n=y\}.
$$
The length distance between $x$ and $y$ is defined by
$$d_\ell(x,y)=\displaystyle\inf_\gamma\{\ell(\gamma):\gamma(0)=x,
\gamma(1)=y\}.
$$

\bigskip

We say $(X,d)$ is a geodesic metric space if for any two points
$x,y\in X$, there exists a curve $\gamma$ connecting $x$ and $y$
such that $\ell(\gamma)=d_\ell(x,y)$. In this case, we say $\gamma$
is a geodesic connecting $x$ and $y$.}

\bigskip

It is clear that a $n$-dimensional complete Riemannian manifold $M^n$ is
always is a geodesic metric space with respect to the metric induced
by the Riemann metric on $M^n$.

Let $(X,d)$ be a geodesic metric space. By a geodesic triangle
$\sigma$ in $X$, we mean a collection of three points $x,y,z$
connected by three geodesics $\sigma_1,\sigma_2,\sigma_3$ and we
call $\sigma_i$ $(i=1,2,3)$ are the sides of geodesic triangle
$\sigma$.

There are several equivalent definitions of Gromov-hyperbolic spaces and
Gromov-hyperbolic groups. We will use the following definition in terms of thin-triangles.

\bigskip

\noindent
\textbf{Definition 3.4.} (Gromov-hyperbolic spaces and groups [Sho91])
{\itshape
(1) Let $(X,d)$ be a geodesic metric space, $(X,d)$ is said to be a
Gromov-hyperbolic space if all geodesic triangles in $(X, d)$ are $\delta$-thin. More precisely,
there is $\delta > 0$ such that each side of any given geodesic triangle in $(X, d)$ is
contained in the $\delta$-neighborhood of the union of the two other
sides. \\

(2) A finitely generated group $\Gamma$ is said to be Gromov-hyperbolic if a Cayley graph
$G_\Gamma$ of $\Gamma$ is Gromov-hyperbolic as a geodesic metric space.}

\bigskip

The relations between a Gromov-hyperbolic space $\tilde M^n$  and its co-compact lattice group $\Gamma$ can be seen by the following proposition.

\bigskip

\noindent
\textbf{Proposition 3.5.} ([Sho91])
{\itshape
Gromov-hyperbolicity is preserved under quasi-isometries. More precisely, suppose that there exists map $f:(X, d) \rightarrow (Y, d_1) $ such that
$\frac{1}{\lambda}d(x,y)-C\leq d_1(f(x),f(y))\leq \lambda d(x,y)+C$
for some constant $C,\lambda>0$ and $f(X)$ is $C$-dense in $Y$. Then
$(X, d)$ is Gromov-hyperbolic  if and only if $ (Y, d_1)$ is Gromov-hyperbolic.}

\bigskip

Suppose that $M^n$ is a compact Riemannian manifold with its fundamental group $\Gamma$.
Then the universal cover $\tilde M^n$ of $M^n$ must be quasi-isometric to any Cayley-graph $G_\Gamma$ of
the fundamental group $\Gamma$.

\bigskip

\noindent
\textbf{Proposition 3.6.}  {\itshape  Let $M^n$ be a compact Riemannian manifold without boundary. Suppose that
$\tilde M^n$ is the universal cover of $M^n$ with the lifted metric and $\Gamma = \pi_1(M^n)$ is the fundamental group
of $M^n$. Then $\Gamma$ is Gromov-hyperbolic if and only if $\tilde M^n$ is a Gromov-hyperbolic metric space.
}

\bigskip

 Under the assumption of nonpositive sectional curvature, our following result shows that Gromov-hyperbolicity condition is
 equivalent to the visibility condition.

 \bigskip
\noindent \textbf{Proposition 3.7.} ([Ca95], [Ca00]) {\itshape Let $M^n$ be a closed  Riemann manifold with nonpositive curvature.
Then $M^n$ is a uniform visibility manifold if and only if its fundamental group $\pi_1(M^n)$ is Gromov-hyperbolic.}

\noindent
\textit{Proof.} For the convenience of readers, we include a short proof here.

\textbf{Step 1.} \textit{To verify that Gromov-hyperbolicity implies visibility.  }

Suppose that $\tilde M^n$ is Gromov-hyperbolic. By Definition 1.1, it is sufficient to verify that ``the universal cover $\tilde
M^n$ contains no totally geodesic embedded Euclidean plane
$\mathbb{R}^2$". Suppose contrary, $\tilde
M^n$ had a totally geodesic embedded Euclidean plane $\Sigma^2
 = \mathbb{R}^2$. We could consider equilateral geodesic triangles $\Delta_\ell$ of side length $\ell$ in $\Sigma^2$.
 As $\ell \to \infty$, the family of triangles $\{ \Delta_\ell \}$ can not be $\delta$-thin for $\ell > 4 \delta$.
 Thus, the geodesic triangles in $\tilde M^n$ can not be uniformly $\delta$-thin for any finite number $\delta$. Hence,
 $\tilde M^n$ can not be a Gromov-hyperbolic space, a contradiction.

\textbf{Step 2.} \textit{To verify that uniform visibility implies Gromov-hyperbolicity.  }
\hfill

We now use an equivalent definition for our uniform visibility
manifold $M^n$. Let us choose $\delta = R(p, \frac{\pi}{2}) =
R(\frac{\pi}{2})$ given by Definition 3.1 for the uniform visibility
manifold $\tilde M^n$. Suppose that $\Delta_{ABC}$ is a geodesic
triangles in $\tilde M^n$ with three vertices $\{A, B, C\}$ and
three sides $\{\sigma_a, \sigma_b, \sigma_c\}$. We are going to show
that our geodesic triangle $\Delta_{ABC}$ is indeed $\delta$-thin
with $\delta = R(\frac{\pi}{2})$. To see this, we may assume that
the side $\sigma_a$ of $\Delta$ has two endpoints $B$ and $C$. For
any interior point point $x \in \sigma_a$ of the side $\sigma_a$, we
draw a new geodesic segment from $x$ to $A$, say $\sigma_{xA}$. It
is clear that
$$
\angle_x(B, A) + \angle_x(A, C) \ge \pi. \eqno(3.6)
$$
Thus, we have
$$
\max\{ \angle_x(B, A), \angle_x(A, C)\} \ge  \frac{\pi}{2}. \eqno(3.7)
$$

  Let $\sigma_{pq}: [0, d(p, q)] \to \tilde M^n$ be the unique geodesic segment from $p$ to $q$ of unit speed. By (3.7), we may assume that
  $$
  \angle_x(B, A) \ge \frac{\pi}{2} \eqno(3.8)
  $$
after re-indexing if needed. We now apply Definition 3.1 to the new geodesic triangle $\Delta_{xAB}$. By (3.8) and our
assumption on the  uniform visibility manifold
$\tilde M^n$, we see that
$$
   d(x, \sigma_{AB}) \le R(\frac{\pi}{2}). \eqno(3.9)
$$

Since $\sigma_c = \sigma_{AB}$, we have verified that ``any point $x$ on one side $\sigma_a$ of the geodesic triangle $\Delta$ is
contained in the $\delta$-neighborhood of the union of the two other
sides $\{ \sigma_c, \sigma_b \}$ for some $\delta = R(\frac{\pi}{2} ) >0$. Therefore, any geodesic triangle $\Delta$ in $\tilde M^n$ is $R(\frac{\pi}{2})$-thin.
It follows from Definition 3.3 that $\tilde M^n$ is Gromov-hyperbolic.

This completes the proof of Proposition 3.7. \qed

\bigskip

 The definition of singular spaces with curvature bounded from above
  $curv \le k$ can be found in
  many textbooks, e.g. [BuBI01] Chapter 9. If a complete metric space $X$ has curvature bounded from above $curv \le k$, then $X$ is called a
  $CAT(k)$-space.

 \bigskip
\noindent \textbf{Proposition 3.8.} ([Bri95]) {\itshape
Let $X$ be a simply-connected and complete metric space of non-positive curvature. Suppose that
there does not exist an isometric embedding of the Euclidean plane into $X$; $X$ has no boundary and suppose that $X$ has a co-compact  lattice
$\Gamma$ such that $X/ \Gamma$ has co-homological dimension $n$ and $H_n(X/\Gamma, \mathbb Z) \neq 0$. Then $\Gamma$ is Gromov-hyperbolic.
 }
 \begin{proof} Proposition 3.7 was extended to possibly singular spaces by Martin Bridson [Bri95]. In fact, the proof of Proposition 3.7 is applicable
 to singular spaces as well.
 \end{proof}

We will discuss the Gromov-norm and minimal volume of visibility manifolds in next section.

\section*{4. Gromov's simplicial norm of visibility n-manifolds and compact 3-manifolds.}
\noindent Our proof of main theorem and theorem 1.3 relies on a
theorem of Igor Mineyev:

 \bigskip
\noindent
\textbf{Theorem 4.1.}
{\itshape (Mineyev [M01])
If $G$ is a hyperbolic group, then the map
$H_b^n(G,\mathbb{R})\rightarrow H^n(G,\mathbb{R})$, induced by
inclusion, must be surjective for $n\geq2$.}

\bigskip

Now we are in position to prove our Main Theorem and Theorem 1.3.

 \bigskip
\noindent \textbf{Proof of Main Theorem:}
\smallskip

Suppose that $M^n$ is a closed visibility manifold. It follows from Proposition 3.7 that the fundamental
group $\pi_1(M^n)$ of $M^n$ is
hyperbolic. It follows from Cartan-Hadamard theorem ([BGS85]) that $M^n$
 must be a closed aspherical manifold ($(K(\pi, 1)$-space).  For any
compact aspherical manifold $M^n$ without boundary, its cohomology ring $H^*(M^n, \mathbb{Z})$ is uniquely determined by its fundamental group
$\Gamma = \pi_1(M^n)$, see Brown's book [Bro94]. In particular, we have
$$
H^*(M^n, \mathbb{Z}) = H^*(\Gamma, \mathbb{Z}), \eqno(4.1)
$$
where $H^*(\Gamma, \mathbb{Z})$ is the co-homology of the group $\Gamma$.

 For Applying Theorem 4.1 to
$\Gamma=\pi_1(M^n)$, we obtain
$H_b^n(\Gamma,\mathbb{R})\rightarrow H^n(\Gamma,\mathbb{R})$ is
surjective. However, $H^n(\Gamma, \mathbb{R})\simeq
H^n(M^n,\mathbb{R})\simeq \mathbb{R}$. It follows that
$H^n_b(\Gamma,\mathbb{R})$ must be nonzero. Thus
$H^n_b(M^n,\mathbb{R})$ is nonzero.  By Gromov's observation (cf.
Proposition 2.8 above), we see that the simplicial volume
$\|M^n\|$ must be non-zero. Main Theorem  now follows from
Gromov's estimate described in Theorem 2.3 above. \qed

\bigskip
\noindent
\textbf{Proof of Main Corollary:} The assertion that the simplicial norm  $\| X/\Gamma \| > 0$ follows
from Proposition 3.8 and Theorem 4.1.   For smooth Riemannian manifolds with Ricci curvature $\ge -(n-1)$, Gromov derived inequality (1.2).
For singular spaces,
Yamaguchi [Yama97] consider the case of curvature $\ge -1$. He further proved that the inequality (1.3) holds.  \qed
\smallskip

\bigskip

We now turn our attention to compact 3-manifolds. In the 3-dimensional case, we can compute the Gromov-simplicial norm for any compact 3-manidolds
due to Perelman's work on Thurston's Geometrization Conjecture. Thurston's proposed decomposition of a 3-manifold has the property that
each summand $N^3_i$ is a locally homogeneous space with possible incompressible
surface boundary. The definition of incompressible surfaces can be recalled
as following.

 \bigskip
\noindent
\textbf{Definition 4.2.}
{\itshape (1) An embedded  two-dimensional sphere $S^2$ in $M^3$ is said to be incompressible if $S^2$ does not bound a $3$-dimensional ball
$B^3$ in $M^3$.

(2) An embedded  two-dimensional torus  $T^2$ in $M^3$ is said to be incompressible
if the inclusion map induces an injective homomorphism from the fundamental group of $T^2$ to the fundamental group of $M^3$,
i.e., the homomorphism $i_*: \pi_1(T^2) \to \pi_1(M^n)$ is injective.

}

\bigskip

There are exactly 8 homogeneous spaces in dimension 3. The list of 3-dimensional
homogeneous spaces includes 8 geometries: $\mathbb{R}^3$, $\mathbb{H}^3$, $\mathbf{S}^3$, $\mathbb{H}^2 \times \mathbb{R}$,
$\mathbf{S}^2 \times \mathbb{R}$, $\tilde{SL}(2, \mathbb{R})$, $Nil$ and $Sol$.

 \bigskip
\noindent
\textbf{Theorem 4.3.} (Perelman [Per02], [Per03])
{\itshape Let $M^3$ be a closed, oriented and  connected
3-manifold. There is  an embedding of a disjoint union $\bigcup \Sigma^2_j $ of incompressible surfaces $\{S^2, \mathbb{R}P^2, T^2, T^2/\mathbb{Z}_2  \}$
 such that
each component of the decomposition $[M^3 -  \bigcup \Sigma^2_j]$ admits a locally homogeneous Riemannian metric of finite volume.  }

\bigskip

It is known that the compact quotient $\mathbb{H}^3/\Gamma$ of 3-dimensional real hyperbolic space $\mathbb{H}^3$ has non-zero Gromov simplicial
norm $\| \mathbb{H}^3/\Gamma \| > 0$, see [Gr82]. Any compact quotient of the remaining seven 3-dimensional homogeneous spaces  $\{
\mathbb{R}^3, \mathbf{S}^3, \mathbb{H}^2 \times \mathbb{R},
\mathbf{S}^2 \times \mathbb{R}, \tilde{SL}(2, \mathbb{R}), Nil,  Sol \}$ is collapsible while keeping sectional curvature bounded from below.
In fact, collapsible 3-manifolds are related to Seifert fiber space and graph-manifolds.

 \bigskip
\noindent \textbf{Definition 4.4.} (Seifert fibered spaces and
graph-manifolds) {\itshape (1) A Seifert fiberation structure on a
compact $3$-manifold $M^3$ is a locally-free circle action on a
finite normal covering $\hat{M}^3$ of $M^3$ such that, denote the covering transformation on $\hat{M}^3$ by $\tau$, we have $\tau(x, \zeta)
= \bar{\zeta} \cdot x $ for all $x \in \hat{M}^3$ and $\zeta \in S^1$. \\

\smallskip
(2) A graph-manifold is a compact $3$-manifold that is connected sum of manifolds each of which is either diffeomorphic to a solid torus or can
be cut apart along a finite collection of incompressible tori into Seifert fibered 3-manifolds.}

\bigskip

The celebrated  Cheeger-Gromov collapsing theory implies that each graph-manifold $M^3$ admits a polarized $F$-structure, cf. [ChG86] and [ChG90].
Hence, for any graph-manifold $M^3$, both minimal volume and simplicial norm of $M^3$ vanish, i.e.,
$$MinVol(M^3) = \| M^3\| = 0. \eqno(4.2)$$

Combining  Perelman's work (Theorem 4.3) and Cheeger-Gromov's collapsing theory, we have

 \bigskip
\noindent \textbf{Corollary 4.5.} (Perelman [Per02], [Per03])
{\itshape  Let $M^3$ be an oriented connected compact 3-manifold
with a Perelman-Thurston the decomposition given by Theorem 4.3. The
simplicial norm of $M^3$ is non-zero if and only if at least one
 of its components of $[M^3 -  \bigcup \Sigma^2_j] $ has hyperbolic geometry. }

\noindent \textit{Proof.} A theorem of Gromov [Gr82] implies that $\| M^3_1 \#_{\cup \Sigma^2_j} M^3_2 \| = \| M^3_1 \| + \| M^3_2\|$, where
$\Sigma^2_j$ are homeomorphic to quotients of $S^2$ or $T^2$. It follows that if $[M^3 - \bigcup \Sigma^2_j]  = N^3_1 \cup N^3_2 ...
\cup N^3_{m}$, then
$$
2 \| M^3 \| = \| \hat{N}^3_1 \| + ... +  \| \hat{N}^3_m\|, \eqno(4.3)
$$
where $\hat{N^3_j}$ is a closed manifold obtained by gluing two copies of $N^3_j$ along their boundaries.

If $\hat{N^3_j}$ is a compact graph-manifold, then $\|\hat{N}^3_j \| = 0$ by a theorem of Cheeger-Gromov [ChG86-90]. If
$\hat{N^3_j}$ is a compact quotient of the remaining seven geometries $\{
\mathbb{R}^3, \mathbf{S}^3, \mathbb{H}^2 \times \mathbb{R},
\mathbf{S}^2 \times \mathbb{R}, \tilde{SL}(2, \mathbb{R}), Nil,  Sol \}$, then $\hat{N}^3_j$ is a graph-manifold. Thus,
$\| \hat{N}^3_j \| \neq 0$ if and only if $N^3_j = \mathbb{H}^3/\Gamma_j$ is diffeomorphic a
quotient of real hyperbolic space $\mathbb{H}^3$ with finite
volume. \hfill \qed

\bigskip

We enclose this paper by the proof of Theorem 1.3.

 \bigskip
\noindent
\textbf{Proof of Theorem 1.3 due to Perelman:}
\smallskip

\noindent
Recall that our 3-manifold $M^3$ has a decomposition: $[M^3 - \bigcup \Sigma^2_j]  = N^3_1 \cup N^3_2 ... \cup N^3_{m}$,  where $\{ \Sigma^2_j\}
$ are either incompressible $2$-spheres $S^2$ or incompressible $2$-tori $T^2$. Since $M^3$ is aspherical, there is no incompressible
$2$-spheres in $M^3$. By our assumption, the fundamental group $\Gamma$ of $M^n$ contains no $\mathbb{Z} \oplus \mathbb{Z}$. Thus, there is no
incompressible $2$-tori in $M^3$ either. Thus, $M^3$ is a locally homogeneous space. Because $M^3$ is aspherical, $M^3$ can not be covered by
$S^3$ or $S^2 \times S^1$. It is known that $Sol$ has no co-compact lattice. If $M^3$ is a compact quotient of $\mathbb{H}^2 \times \mathbb{R}$, then, by
a theorem of Eberlein [Eb83], a finite normal cover of $M^3$ is diffeomorphic to $\Sigma^2 \times S^1$,
where $\Sigma$ is a closed surface of genus $\ge 1$. In summary, we have the following fact.

\smallskip
\noindent
\textbf{Fact 4.6.} If $M^3$ is a compact quotient of $\{
\mathbb{R}^3,  \mathbb{H}^2 \times \mathbb{R},  \tilde{SL}(2, \mathbb{R}), Nil \}$, then the fundamental group
$\Gamma $  of $M^3$ has a subgroup $\hat\Gamma$
of finite index such that $\hat\Gamma$ has a nontrivial center $\mathcal C$ containing $\mathbb{Z}$.
\smallskip

The verification of Fact 4.6 goes as follows. (i) For $\mathbb{R}^3$, the Biebarbach theorem (cf. [BaGS85]) states that any co-compact lattice
$\Gamma$ contains a subgroup $\hat\Gamma$ isomorphic to $\mathbb{Z} \oplus \mathbb{Z} \oplus \mathbb{Z}$. (ii) We already discussed the compact
quotient of $\mathbb{H}^2 \times \mathbb{R}$ by using Eberlein's result. (iii) It is known that $\tilde{SL}(2, \mathbb{R})$ is isometric to
the unit tangent bundle $S\mathbb{H}^2$ of the Poincare disk. There is a fiberation $ S^1 \to S\mathbb{H}^2 \to \mathbb{H}^2$. The fundamental subgroup
corresponding to the $S^1$ is the center of $\pi_1(S\mathbb{H}^2 /\Gamma     )$. For the nilpotent group $Nil$, it is well-known that the co-compact
group
$\pi_1( Nil/\Gamma)$ has a subgroup $\hat\Gamma$
of finite index such that $\hat\Gamma$ has a nontrivial center isomorphic to $\mathbb{Z}$.

Since the cohomological dimension of $M^3$ is equal to $3$ (cf. [Bro94]), its fundamental group $\Gamma$ can not be isomorphic to $\mathbb{Z}$.
Thus, by Fact 4.6, we see that if $M^3$ is a compact quotient of $\{
\mathbb{R}^3,  \mathbb{H}^2 \times \mathbb{R},  \tilde{SL}(2, \mathbb{R}), Nil \}$, then its fundamental group $\Gamma$ contains a subgroup
isomorphic to $\mathbb{Z}\oplus \mathbb{Z}$, which contradicts to our assumption.

Therefore, by our assumptions on $M^3$, there is  only  one possibility: $M^3$ is diffeomorphic to a compact quotient of $\mathbb{H}^3$.
Thus, $MinVol(M^3) \ge \frac{1}{24} \|M^3\| > 0$.

Recall that, by a theorem of Yau [Y71], the fundamental group of a visibility manifold contains no subgroup isomorphic to $\mathbb{Z}\oplus \mathbb{Z}$.
The above argument shows that $M^3$ is diffeomorphic to a compact quotient $\mathbb{H}^3/\Gamma$ of $\mathbb{H}^3$.

This completes the proof of Theorem 1.3. \hfill \qed

 \bigskip
\noindent
\textbf{Corollary 4.7.}
{\itshape Let $M^3$ be a closed $3$-dimensional aspherical manifold. Suppose that the fundamental group $\pi_1(M^3)$ contains
no sub-group isomorphic to $\mathbb{Z} \oplus \mathbb{Z}$, suppose that $g$ is a smooth metric of scalar curvature $\ge -6$.
Then $M^3$ is diffeomorphic to a compact quotient of $\mathbb{H}^3$ and the volume of $(M^3, g)$ has a lower bound:
$$
Vol(M^3, g) \ge Vol(\mathbb{H}^3/\Gamma) = c_3 \| \mathbb{H}^3/\Gamma \| > 0.
$$
}

  Corollary 4.7 is a direct consequence of Theorem 1.3 and Theorem 2.5.

\vspace{5mm}
\noindent
{\bf Acknowledgement:} The second author would like to express his gratitude to his advisor ---- Professor Xuan-Guo Huang for teaching.
Authors want to thank Professor Viktor Schroeder for his summer course on ``Asymptotic Geometry" at Nanjing University,
China in 2008.
 We are also grateful to Nanjing University for its hospitality. Professor D. Kotschick kindly bring the work of Besson-Courtois-Gallot
 [BCG96] to our attention, which lead to Corollary 4.7.
 Professor Frank Morgan and Professor Gudlaugur Thorbergsson made a number of suggestions to improve the exposition of our paper.

    We are very grateful to Igor Belegradek for bring the references [Bri95] and [Yama97] to our attention. In 1990-1991, among other things,
    the first author
    lectured Proposition 1.14 of [Ca95] in a graduate topic course at Cornell University, for which we felt very flattered that Martin Bridson
     was attending. We
    are also very delighted to find out that M. Bridson was able to extend Proposition 1.14 of [Ca95] on smooth manifolds to
    a result on singular spaces in [Bri95] after his graduation from Cornell.

 Finally, we indebted to
 Professor Francois Ledrappier and Professor Karsten Grove for several conservations on Gromov's simplicial volumes.

\bigskip

\noindent
{\small Jianguo Cao }\\
{\small Mathematics Department} \\
{\small University of Notre Dame}\\
{\small Notre Dame, IN  46556, USA} \\
{\small cao.7@nd.edu} \\
\\
{\small Mathematics Department} \\
{\small Nanjing University}\\
{\small Nanjing, 210093, China} \\
\\
Xiaoyang Chen \\
{\small School of Mathematical Sciences}\\
{\small Fudan University} \\
{\small Shanghai, 200433, China}\\
{\small 071018012@fudan.edu.cn}

\end{document}